\documentclass[aap,preprint]{imsart}

\RequirePackage[OT1]{fontenc}
\RequirePackage{amsthm,amsmath}
\RequirePackage[numbers]{natbib}
\RequirePackage[colorlinks,citecolor=blue,urlcolor=blue]{hyperref}

\arxiv{arXiv:1410.4814}

\startlocaldefs
\numberwithin{equation}{section}
\theoremstyle{plain}

\endlocaldefs

\newtheorem{theorem}{Theorem}[section]
\newtheorem{lemma}[theorem]{Lemma}
\newtheorem{proposition}[theorem]{Proposition}

\newtheorem{definition}[theorem]{Definition}
\newtheorem{corollary}[theorem]{Corollary}
\newtheorem{remark}[theorem]{Remark}

\newtheorem{hyp}[theorem]{Hypothesis}

\newcommand{\cB}{{\cal B}}

\newcommand{\cC}{{\cal C}}

\newcommand{\cL}{{\cal L}}
\newcommand{\cX}{{\cal X}}

\def \a {{\alpha}}

\def \e {{\varepsilon}}

\def \m {{\mu}}

\def \s {{\sigma}}

\def \t {{\tau}}

\def \d {{\delta}}
\def \p {{\pi}}

\newcommand{\be}[1]{\begin{equation}\label{#1}}
\newcommand{\ee}{\end{equation}}

\newcommand{\bl}[1]{\begin{lemma}\label{#1}}
\newcommand{\el}{\end{lemma}}

\newcommand{\ind}{\mathbf{1}}

\newcommand{\br}[1]{\begin{remark}\label{#1}}
\newcommand{\er}{\end{remark}}

\newcommand{\bt}[1]{\begin{theorem}\label{#1}}
\newcommand{\et}{\end{theorem}}

\newcommand{\bd}[1]{\begin{definition}\label{#1}}
\newcommand{\ed}{\end{definition}}

\newcommand{\bcl}[1]{\begin{claim}\label{#1}}
\newcommand{\ecl}{\end{claim}}

\newcommand{\bp}[1]{\begin{proposition}\label{#1}}
\newcommand{\ep}{\end{proposition}}

\newcommand{\bc}[1]{\begin{corollary}\label{#1}}
\newcommand{\ec}{\end{corollary}}

\newcommand{\bi}{\begin{itemize}}
\newcommand{\ei}{\end{itemize}}

\newcommand{\ben}{\begin{enumerate}}
\newcommand{\een}{\end{enumerate}}

\begin {document}

\begin{frontmatter}
\title{Conditioned, quasi-stationary, restricted measures and escape from metastable states}
\runtitle{Conditioned, quasi-stationary, restricted measures and escape from metastable states}

\begin{aug}
\author{\fnms{R.} \snm{Fernandez}\ead[label=e1]{R. Fernandez1@uu. nl}},
\author{\fnms{F.} \snm{Manzo}\ead[label=e2]{manzo.fra@gmail.com}}
\author{\fnms{F.} \snm{Nardi}\ead[label=e3]{f.r.nardi@tue.nl}},
\author{\fnms{E.} \snm{Scoppola}\ead[label=e4]{scoppola@mat.uniroma3.it}}
\and
\author{\fnms{J.} \snm{Sohier}\ead[label=e5]{jsohier@tue.nl}}
%
\runauthor{R. Fernandez et al.}


\address{ Roberto Fernandez \\
University of Utrecht \\ 
 PO Box 80010 \\
 3508 TA Utrecht \\
 The Netherlands \\
\printead{e1}\\}

\address{ Francesco Manzo \\
Department of Mathematics\\
University of Roma Tre\\
Largo San Murialdo 1\\
00146 Roma \\
Italy\\
\printead{e2}\\}

\address{ Francesca Nardi \\
Eindhoven University of Technology \\
Department of Mathematics and Computer Science \\
MetaForum (MF) 4.076 \\
P.O. Box 513 \\
5600 MB Eindhoven \\
The Netherlands \\
\printead{e3}\\}

\address{ Elisabetta Scoppola \\
Department of Mathematics\\
University of Roma Tre\\
Largo San Murialdo 1\\
00146 Roma \\
Italy\\
\printead{e4}\\}

\address{ Julien Sohier \\
Eindhoven University of Technology \\
Department of Mathematics and Computer Science \\
MetaForum (MF) 4.118 \\
P.O. Box 513 \\
5600 MB Eindhoven \\
The Netherlands \\
\printead{e5}\\}

\end{aug}

\begin{abstract}
  We study the asymptotic hitting time $\tau^{(n)}$ of a family of Markov  processes $X^{(n)}$ to
  a target set $G^{(n)}$ when the process starts from a ¨trap¨ defined by very general properties.
  We give an explicit description of the law of $X^{(n)}$ conditioned to stay within the
  trap, and from this we  deduce the exponential distribution of $\tau^{(n)}$. Our approach is
  very broad ---it does not require reversibility, the target $G$ does not need to be a rare event, and 
  the traps and the limit on $n$ can be of very general nature--- and leads to explicit bounds on
  the deviations of $\tau^{(n)}$ from exponentially.   We provide two non trivial 
      examples to which our techniques directly apply. 
  
\end{abstract}

\begin{keyword}[class=MSC]
\kwd{60J27}
\kwd{60J28}
\kwd{82C05.}
\end{keyword}

\begin{keyword}
\kwd{ Metastability, continuous time Markov chains on discrete spaces, hitting times, asymptotic exponential behavior.}
\end{keyword}

\end{frontmatter}

%

%
%
%
%
%
%
%
%
%
%
%
%

\section{Introduction}

\subsection{Scope of the paper}

Metastability and related phenomena are  associated to systems "trapped" for a long time in some part of their phase space, from which they
emerge in a random time ---the \emph{exit time}--- expected to have an asymptotic exponential distribution. 
They are the subject of many current studies in the mathematical and mathematical physics literature.
This recent activity is motivated, on the one hand, by the confirmation of their widespread occurrence in most evolutions and,
on the other hand, on the emergence of a host of new techniques for their rigorous analysis ---cycle decomposition
\cite{CaCe,Ca,CGOV,OS1,OS2,OV,MNOS,Trouve1},  potential theoretic techniques \cite{BovMetaPTA,BEGK00,BEGK02,BG},
renormalization \cite{Scop93,Scop94}, martingale approach \cite{BLM,BL1,BL2}.  
 In this article, we focus on an essential component
of a metastability analysis:  the description of the state of a system trapped in a region $A$ and the estimation of the law of the exit 
time from this (meta)stable trap to a target set $G$.  These times are the
building blocks for studies of evolutions through arbitrarily complex scenarios involving many of these traps.  Our treatment 
is inspired by the first rigorous paradigm proposed for this type of studies ---the  \emph{path-wise approach}  introduced
in  \cite{CGOV} and developed, for instance, in \cite{OS1,OS2,OV,MNOS}. 

In the remaining part  of the introduction we discuss in detail the main advantages of our approach, but here is a brief summary:
(i) We do not make any assumption on the nature of the limit process involved, in particular it applies to fixed-temperature, 
infinite-volume limits. (ii) The process is not assumed to be reversible, and the traps do not need
to be narrowly supported.
(iii) (Meta)stable states are described by measures, not just through a few reference configurations.
Furthermore, the  different natural candidate measures are shown to provide equivalent descriptions within explicit error bounds.
(iv) There is no constraint on the exterior $G$ of the trap ---in particular, it does neither need to be a rare set, nor to be associated 
 to a stable state.
(v) Exit laws are precisely estimated with explicit bounds on deviations from exponentiality. 
(vi) Our approach relies on a novel type of proof based on controlling the proximity to the quasistationary measure.  This
simplifies proofs and strengthens results.

\subsection{Issues addressed by this approach}\label{notiss}

Let us first discuss in some detail the issues we deal with.

\paragraph{General traps} In many visualisations of a trap, authors have in mind some energy
profile associated to the invariant measure of the process defining
the evolution.  A trap corresponds, in this case, to an energy well, and metastability refers to the exit from such a well to 
a deeper well, leading to stability, that is, corresponding to the asymptotic support of the invariant measure. 
This point of view is fully exploited, for instance, in the
reference book \cite{OV}. It is important to keep in mind, however, that this is only one of many possible metastability
(or tunnelling) scenarios.  Indeed, in many simple examples the invariant measure is not linked to any useful ``energy" profile.
For instance, any shuffle algorithm for a pack of cards leaves the uniform measure invariant.  This measure corresponds to 
"flat", well-less energy profile.  Yet, in this paper we show that for the 
 well known Top-In-At-Random model, the time it takes an initially well shuffled pack to attain a particular order
is exponentially 
distributed (see section \ref{sec3:2}).  In this example the "well" $A$ is entropic in nature.  More generally, processes can
define free-energy wells. We stress that our setup relies on hypotheses insensitive to the nature of the "well".

\paragraph{Measures associated to traps}
States are probability measures and hence transitions between trapped states have to correspond to transitions between measures (asymptotically)
supported  by these traps.  Nevertheless, 
in many instances metastable states are associated to individual configurations marking the ``bottom of the trap".  Such studies
are suited to traps that become, asymptotically, abrupt enough to allow only small fluctuations in the trapped state.  A more general
scenario should involve "wide traps" corresponding to measures with an extended support and a corresponding thermalization scale.
 For a given trap, there are three natural candidate measures: the restriction of
the invariant measure, the quasistationary measure and the empirical measure of the process before exiting the trap. 
These measures have complementary properties; in particular the quasistationary measure automatically leads to exponential exit times.  
In this paper, we shed light on the links between these measures in our framework. 
 
 \paragraph{General assumptions}
Our approach for metastability includes the following features. 
\begin{itemize}
\item \emph{No need of reversibility:} Many approaches depend crucially of the reversible character of the invariant measure.
This is, in particular,
true in the funding literature of the subject by Keilson \cite{Ke} and Aldous and Brown  \cite{AB,AB2}. Until the recent paper \cite{La}, where potential
 theoretic tools were used to prove the metastable behavior for a non-reversible dynamics, it was also a 
prerequisite for the application
of potential theory.  Nevertheless,
many examples (we deal with one such example in Section \ref{sec3:2}) show that this hypothesis is not necessary to
observe exponential escape times.

\item\emph{General nature of the target set $G$:}
Exit times are defined by hitting times to a certain "boundary set" $G$ which, in many instances, is associated to the 
\emph{saddle points}, or bottlenecks that the system has to cross on its exit trajectories.  As such, it is often assumed 
that the set $G$ is asymptotically negligible.  Alternative metastability (and tunnelling) studies use $G$ as the support 
of the stable measure (=bottom of deeper wells).  In these cases, the set $G$ has asymptotically full measure, 
 which is in sharp contrast with the previous ones.  These two 
cases show that, for general exit times studies, the set $G$ should simply be associated with the exterior of the well,
without any assumption on its asymptotic measure.
\end{itemize}
\medskip

\paragraph{The nature of the asymptotic regime}
Metastability, and the exponential escape law, only appear asymptotically in appropriate parameters that gauge the 
nature and depth of the trap.
The type of parametrization determines how the metastable regime is approached.  Some approaches ---e.g.\ 
low-temperature limits--- involve traps that 
become asymptotically more abrupt.  Others ---e.g. card shuffling of larger packs--- keep the geometry of 
the trap fixed but make the external set $G$ progressively farther from the "bottom" of the trap.
In more complicated situations the limit involves both a change in depth and the complexity of the traps. This 
happens, for instance, in the study of fixed-temperature spin systems in the thermodynamic limit, for which
the changes in the volume lead to a more complex trap scenario with an associated proliferation ( "entropy effect")
of escape routes (e.g. location of "critical droplets").  
Our approach applies equally well to these different limits.

\paragraph{The type of transition probabilities}
In many standard metastability studies, the stochastic evolution is itself defined through an energy function (identical or closely 
related to the one defined by the invariant measure).  This is the case, for instance, in the Metropolis algorithm, in which
transition probabilities are proportional to the exponential of the positive part of energy differences.  A related
feature of these models is the
fact that the "barriers" between traps ---determining, for instance, the mean escape time--- are proportional to the energy differences
between configurations marking the bottom and the top of energy wells. In the framework of statistical mechanics models, 
an important family of models  
 excluded by this setup are the 
 cellular automata \cite{CN03}. It turns out that the parallel character of their dynamics leads to the existence
 of a large number of paths between any
two configurations.  The probabilistic barrier characterizing each of these paths cannot be described only by energy differences. 

 The particular form of the transition probabilities plays no special role in our general approach.

\subsection{Main features of our approach}

Here is how we deal with the issues discussed above.

\paragraph{General traps and asymptotic regime}
We consider a family of continuous time irreducible Markov chains $X^{(n)}_t$ on finite
state spaces $\cX^{(n)}$.  The asymptotic regime is associated to $n\to\infty$,  but we do not need to specify the nature of
the parameters involved. Particular relevant examples are the case where the limit may involve the divergence of the 
cardinality of the state spaces, and/or
some parameter in the transition probabilities such as the inverse temperature.  For each $n$ the state space is divided
into a \emph{trap} $A^{(n)}$ and its complementary set $G^{(n)}=\cC^{(n)}\setminus A^{(n)}$.  We do not assume any particular property 
of $G^{(n)}$; exit times are only determined by the evolution inside $A^{(n)}$ and the structure of $G^{(n)}$ is 
irrelevant for the exit times.  The traps include sets $B^{(n)}\subset A^{(n)}$ of configurations ---associated to ``bottoms 
of the well"--- which can have arbitrary size as long as they satisfy certain natural assumptions.

More precisely, the definition of trap is contained in three physically natural hypotheses ---spelled out in Section \ref{sec2.1}
below--- that can be roughly described as follows:

\begin{description}
\item{\emph{Fast recurrence:}} For any initial configuration, the process falls within a controlled time $R^{(n)}$ either in $B^{(n)}$ or 
in $G^{(n)}$. This \emph{recurrence time} acts as a reference time for the whole of the analysis. 
\item{\emph{Slow escape:}} Starting from configurations in $B^{(n)}$, the time $\tau^{(n)}$ the process takes to hit $G^{(n)}$ is much larger
than the recurrence time $R^{(n)}$.
\item{\emph{Fast thermalization:}} A process started in $B^{(n)}$  achieves
"local thermalization" within $B^{(n)}$--- in a time much shorter than the escape time $\tau^{(n)}$.
\end{description}
The two time scales ---escape and local thermalization---  actually define the trap and the limit process.  This definition
is sufficiently general to accommodate different types of traps (e.g.\ energy, entropy- or free-energy-driven) and 
asymptotic limits (e.g.\ low temperature, large volume or combinations of both).  The scales do not need
to be asymptotically different, a fact that allows the consideration of less abrupt traps.  
The recurrence hypothesis controls, in particular, the presence of  complicated landscapes  in $A^{(n)}\setminus B^{(n)}$.  As discussed below,
this control is
  not necessary if $G^{(n)}$ is a rare set.  Furthermore, in concrete examples it is usually not
  hard to achieve recurrence at some appropriate time scale, due to the 
  submultiplicative property of the probability of non-recurrence
   \begin{equation*}
      \sup_{x\in\cX^{(n)}}\mathbf{P}(\t^x_{G^{(n)}\cup B^{(n)}}>kt)\le\sup_{x\in\cX^{(n)}}\mathbf{P}(\t^x_{G^{(n)}\cup B^{(n)}}>t)^k.
   \end{equation*}

\paragraph{Metastable states}
Our approach is based on the fact that, with the above definition of a trap, the three possible notions of trapped state ---the
restriction of
the invariant measure, the quasistationary measure and the empirical measure of the process before exiting the trap---  are
largely equivalent and can be used indistinctively.  In fact, our proof provides explicit bounds on the asymptotic behavior
of their distances.  Therefore,  metastable states can be defined by any of them and are endowed with the physical
properties associated to each of them.  Mathematically, this equivalence leads to precise quantitative
 estimates on the distribution of the first hitting time of $G^{(n)}$, through notoriously simplified proofs.
 
 \paragraph{General target set $G^{(n)}$}
Except for the recurrence hypothesis, the structure of the set $G^{(n)}$ plays no role in our theory.  It can, therefore,
correspond to rare or non-rare parts of the configuration space.  Furthermore, this generality makes
our results applicable to a wide type of phenomena.  Indeed, in a truly
metastability scenario the trap corresponds to an asymptotically unlikely family of configurations and the escape is
followed by a fall into the true support of the stable phase.  In other cases, however, the escape may be followed by
a visit to an equally likely
set of configurations; a situation more appropriately referred to as \emph{tunnelling}.  This phenomenon occurs, for instance, if stable
states are supported in several asymptotically disconnected regions and interest focuses on the way the system migrates from one of
these regions to another one. Hence, the insensitivity to the scenario after the exit leads to results  useful for the study of
metastability, tunnelling or any evolution involving traps.  

\paragraph{Application to non-reversible dynamics}
Other than the preceding assumptions, we do not assume any particular property of the transition probabilities, nor do we
assume reversibility.  This work is, hence, part of the so far quite reduced circle of
metastability results in the non-reversible case.  Until the work \cite{La}, lack of reversibility was precluding the potential
theoretical approach, and relevant publications are very recent  \cite{GL,BL2}.  In particular, the 
convergence in law of escape times has  been obtained in  \cite{BLM}, but assuming that $G^{(n)}$ is a rare event.
In \cite{OI}, the author provides a control on tail distributions  that is close to ours, yielding an error of the form 
 $\exp(O(\e)t)$ instead of our $O(\e)\exp(-t)$ (see also (ii) in Remark \ref{remdop}).  

  The technique of cycle decompositions 
\cite{OS2,CaCe} applies also to non-reversible dynamics  for chains with exponentially small
transition probabilities at fixed-volume  in the low-temperature   limits.  

 \paragraph{ General assumptions on the asymptotic regime} We do not make any assumption on the nature of the limiting process 
  involved. In particular, we cover the case of fixed temperature and infinite volume limit, where it naturally arises an 
  important entropic contribution. 
  There are only few studies that deal with the latter for specific models under reversible dynamics. They consider the case 
   in which the volume grows exponentially with the inverse temperature \cite{GHNOS2, BdHS, GN}. In \cite{SchSh}, the authors consider 
    the kinetic Ising model at infinite volume, fixed temperature and asymptotically vanishing external magnetic field.

\bigskip

\paragraph{Summary of results} For systems satisfying the three hypotheses above, we prove that the first hitting times to $G^{(n)}$ 
are asymptotically exponentially distributed, within explicitly bounded error terms.  Furthermore, we get explicit control
of the variational distance between the quasi stationary measure and
 the evolved measure of the process conditioned to stay within $A^{(n)}$.  The exponential
law is subsequently deduced from the properties of the quasi-stationary measure.  This approach leads to
stronger control on the distribution of the
hitting time than the one obtained in \cite{OI}.  In fact, our estimations are valid \emph{for every
starting point} in $A^{(n)}$, and hence yield a precise description
 of the empirical measure of the chain
up to its first escape from $A^{(n)}$.  This approach is different ---and more effective--- than the one adopted
in a precedent paper  \cite{FMNS}, which was restricted to traps
       characterized by a single-configuration bottom.

Intuitively, the picture we obtain is the following:  in a time of the order of the recurrence time $R^{(n)}$ the 
process chooses, depending on its starting point, whether it reaches directly $G^{(n)}$ or 
  whether it visits $B^{(n)}$ before.  Once in $B^{(n)}$, the process reaches in a very short time a metastable equilibrium 
  within $A^{(n)}$ that is accurately described by the quasistationary measure.  Due to Markovianness and the long mean 
  exit times, the process requires an exponential time to exit from this temporary equilibrium.

      In the second part of the paper, we illustrate the power of our results through applications  to two 
      interesting models.  Both models have features that put them outside most general theories. 
      The first model is a  reversible birth and death process issued from the context of dynamical polymer models. 
      It has the particularity that its energy profile exhibits a double well, but separated by an energy barrier
      that is only logarithmic.   Our second example is the well known Top In At Random 
      (TIAR) model, which has a uniform invariant measure.  This is an example of a non reversible dynamics with
      a purely entropic barrier.  As far as we know, this is the first time this model has been shown to exhibit
      exponential behavior of the hitting time to a rare set.  
          
           Outline of the paper: in Section \ref{Not} we give the main definitions and recall
           some well known results, and in Section 
   \ref{Res} we precise the hypotheses and detail our results. Applications are discussed in Section \ref{Ex}. Finally,
   Section \ref{Pro} is devoted to the proofs.

\section{Setting}\label{Not}
 
Let $\cX^{(n)}$ be a sequence of finite state space depending on a parameter $n$ and $X^{(n)}_t$ be a sequence
of {\it continuous  time irreducible  Markov chains} on them. We are interested in the asymptotics  $n\to\infty$. We should 
 typically think of $n$ as being related to: 
 \begin{enumerate}
  \item the cardinality of $|\cX^{(n)}|$. This is typically the case for large size physical systems, and our
  hope would be that the techniques developed in this paper apply to the infinite volume of well known models issued from statistical
  physics. 
  \item the inverse temperature in a Freidlin-Wentzell setup. This issue has been very studied, 
  classical references by now are \cite{OV} and \cite{Ca}. 
 \end{enumerate}

  We denote
by $Q^{(n)}$ the matrix of transition rates on $\cX^{(n)}$ generating the chain $X^{(n)}_t$ and
 by $(\pi^{(n)}, \, n\ge 1)$ the corresponding sequence of  invariant measures. The family of Markov chains $X^{(n)}_t$ 
  could be equivalently defined in terms of their kernels 
   \begin{equation}
    P^{(n)} = \ind + Q^{(n)},
   \end{equation} 
    and via their continuous time semi group 
   \begin{equation}
    H_{t}^{(n)} f := e^{t Q^{(n)}} f = e^{-t} \sum_{k \geq 0} \frac{t^{k} (P^{(n)})^{k}}{k!} f.
   \end{equation}

  We do not assume reversibility so that the  adjoint kernels  $(P^{(n)})^*$ defined by
   \begin{equation}\label{rev}
      (P^{(n)})^*(x,y)={\p^{(n)}(y)\over \p^{(n)}(x)}P^{(n)}(y,x)
   \end{equation} 
  do not necessarily coincide with $ P^{(n)}$. We will denote by $ X^{\leftarrow(n)}_t$ the corresponding
  time reversal processes.

Discrete time Markov chains are covered by our techniques and we could also  consider more general cases but for
an easier presentation of our results, we
prefer to restrict our analysis to this setup.

We denote by $\mathbf{P}(.)$ and $\mathbf{E}(.)$ generic probabilities and mean values. We will specify
in the events and in the random variables the initial conditions of the process, 
so that for instance, for $A^{(n)}\subset \cX^{(n)}$ we define the {\it first hitting time to} $A^{(n)}$
for the chain  $X^{(n)}_t$ starting  $x\in\cX^{(n)}$:
\be{hitting}
\t^{(n),x}_{A^{(n)}}=\inf \{t\ge 0 :\; X^{(n),x}_t\in A^{(n)}\}.
\ee

 We also use the small-o notation $o(1)$ to denote a positive quantity which has limit $0$ as $n\to\infty$.

 We will prove our results keeping in mind the asymptotic regime $n \to \infty$. Nevertheless, we stress that most of our estimates go beyond 
  asymptotics, and more precisely we were able to show our main results with explicit error terms. We choose to write the 
   major part of this article for a fixed value of $n$ (and as a consequence and for lightness of notation, from now on, we
   drop this explicit dependence from the notations, except in some particular cases like Corollary \ref{convlaw}), denoting by small latin letters quantities which should be thought of as 
    asymptotically $o(1)$. For example, the quantity $f$ in \eqref{f} should be thought of as a sequence $f_{n} = o(1)$. 
 
For positive sequences $a_{n}, b_{n}$, we will also use the standard notations $a_{n} \sim b_{n} $  as 
soon as $\lim_{n \to \infty} a_{n}/b_{n} = 1$, and $a_{n} \gg b_{n}$ if $b_{n}/a_{n} = o(1)$.

The notation $X=\left(X_t \right)_{t\in\mathbf{R}}$ is used for  a generic chain on $\cX$.

\subsection{Measures and distances}
\label{diffmeas}
Let $X_t$ be an irreducible continuous time Markov chain on $\cX$ with stationary measure $\pi$
and transition rates $Q$, and
let $A\subset \cX$ be a given set. From now on, we denote by $G$ the complementary of $A$, and we define the following measures:
\bi
\item[-] {\bf The  invariant  measure restricted on $A$:}
\be{misristr}
\pi_A(.):= \frac{1}{\pi(A)}\pi(.)\ind_{\{.\in A\}}.
\ee

\item[-] {\bf The measure of evolution: }
for any $t\ge 0$ and  $x\in\cX$ let 
\be{misevol}
\mu^x_t(.):=\mathbf{P}(X^x_t=. ),
\ee
and more generally for any probability measure $\nu$ on $\cX$: 
\be{evolpa}
\mu^{\nu}_t(.):=\mathbf{P}(X^{\nu}_t=. ):=\sum_{x\in \cX}\nu(x)\mathbf{P}(X^x_t=. ).
\ee
\item[-] {\bf The conditioned evolution measure on $A$:}

for any $t\ge 0$ and $x\in A$ the evolution measure conditioned to $A$ is defined as:
\be{miscond}
\tilde \mu^x_{A,t}(.):=\mathbf{P}(X^x_t=. |\t^x_{A^c}>t)= \frac{\mathbf{P}(X^x_t=. ,\t^x_{A^c}>t)}{\mathbf{P}(\t^x_{A^c}>t)},
\ee
and more generally for any probability measure $\nu$ on $\cX$: 
\be{condpa}
\tilde \mu^{\nu}_{A,t}(.):=\sum_{x\in\cX} \nu(x)\tilde\mu^x_{A,t}(.).
\ee
\item[-] {\bf Quasi stationary measure on $A$:}

 The quasistationary measure $\m^*_A(.)$ is a classical notion; it can be defined for example in the following way: 
 \be{misqs}
\m^*_A(.):=\lim_{t\to \infty}\tilde \mu^{\p_A}_{A,t}(.).
\ee

We recall the following elementary property about $\mu_{A}^{*}$; this is a classical result which can be
 found for example in 
  \cite{Sen}, and which states that starting from  $\mu_{A}^{*}$, the exit law from $A$ is exactly exponential (no correction is needed). Here 
   and in the rest of the paper, we set $ \mathbf{E}[ \tau^{\mu^{*}_{A}}_{G}] =: T^{*}$. 
  
   \begin{proposition}\label{exppro}
    For any $t \geq 0$, the following equality holds:
    \begin{equation}
     \mathbf{P}[ \tau^{\mu_{A}^{*}}_{G}/T^{*} > t ] = e^{-t}.
    \end{equation} 
     
   \end{proposition}
   
   \item[-] {\bf The empirical measure on $A$:}
   
   for any $x\in A$, by denoting with $\xi^x_G(y)$ the local time spent in $y$ before $\t_G$ starting at $x$
   we can  define the empirical measure:
   \be{emme}
   \m^{(em) x}_A(y)=\frac{\mathbf{E}\xi^x_G(y)}{\mathbf{E}\t^x_G}.
   \ee

\item[-] {\bf Total variation distance:} given two probability measures $\pi_{1}$ and $\pi_{2}$ on $\cX$, we define:
$$
d_{TV}(\p_1,\p_2):= {1\over 2}\sum_{x\in\cX}|\p_1(x)-\p_2(x)|=\max_{A \subset \cX}|\p_1(A)-\p_2(A)|,
$$
and  
\be{dt}
d(t)=\max_{x\in\cX}d_{TV}(\mu^x_t,\pi),\qquad \bar d(t):=\max_{x,x'}d_{TV}(\mu^x_t,\mu^{x'}_t).
\ee
 \ei

It is well known (see for instance \cite{AD}) that $\bar d$ is submultiplicative; namely, for $m,n>1$
\be{submult}
\bar d(m+n)\le \bar d(m)\bar d(n), 
\ee
and
$$
\bar d(t)\le d(t)\le 2 \bar d(t).
$$

 For a given set $K \subset \cX$, we introduce the quantity
 \begin{equation}\label{defbardAt}
  \bar d_K(t):=\max_{x,x'\in K}d_{TV}(\mu^x_t,\mu^{x'}_t).
 \end{equation} 
 
  A submultiplicative property such as \eqref{submult} does not hold in general for the quantity $\bar d_K(t)$.

\section{Results}\label{Res}
We present in this section our main results
on the asymptotic vicinity of the conditioned measure and the conditioned invariant one.
  As a corollary we can prove exponential behavior of the first hitting time of $G$ with accurate control on the error terms.

   \subsection{ Hypotheses}\label{sec2.1}

Let  $G\subset \cX$ be a set representing "goals" for the chain $X_t$,
and  define $A:=\cX\backslash G$, so that  $\tau^{x}_{G}$ represents the first exit time from $A$ and
$\tau^{\leftarrow x}_{G}$ is the first exit time from $A$ for the time reversal process $ X^\leftarrow_t$.
To observe asymptotic exponentiality of the hitting time of $G$, we need to make a couple of assumptions on the set $A$; define 
the quantity 

\begin{equation}\label{f1}
 f_A(t):= \mathbf{P}(\tau^{ \pi_{A}}_{G}\le t)={1\over 2}\sum_{x\in A}\pi_A(x)\Big[ \mathbf{P}(\tau^{ x}_{G}\le t)+\mathbf{P}(\tau^{\leftarrow x}_{G}\le t)\Big].
\end{equation} 

Indeed one readily realizes that $\mathbf{P}(\tau^{ \pi_{A}}_{G} > t) = \mathbf{P}(\tau^{\leftarrow  \pi_{A}}_{G} > t) $ using time reversal :
\begin{equation}
 \begin{aligned}
  \mathbf{P}(\tau^{\leftarrow  \pi_{A}}_{G}\le t) & = \sum_{x,y \in A} \sum_{x_{1},\ldots,x_{t-1} \in A} 
  \frac{\pi(x)}{\pi(A)} P^{*}(x,x_{1}) \ldots P^{*}(x_{t-2},x_{t-1}) P^{*}(x_{t-1},y) \\ 
   & = \sum_{y,x \in A} \sum_{x_{1},\ldots,x_{t-1} \in A} 
  \frac{\pi(y)}{\pi(A)} P(y,x_{t-1}) P(x_{t-1},x_{t-2}) \ldots  P(x_{1},x) \\
   & = \mathbf{P}(\tau^{ \pi_{A}}_{G}\le t)
 \end{aligned}
\end{equation} 
 and \eqref{f1} follows. 
 
 Note that we carried this proof in the discrete time setup for lightness of notations, the continuous time case works in the same way.

  We need three natural assumptions on the behavior of the process. 
  \begin{enumerate}
   \item[\textbf{E(R,f)}] \textbf{Slow escape}: on a time scale $2R$, starting from the invariant 
    measure restricted to the set $A$, the process  hits $G$ with small probability; concretely,
    for $f\in (0,1)$  and $R>1$, we say that the set $A$ satisfies the hypothesis $E(R,f)$ if
\begin{equation}\label{f}
  f_A(2R)\le f.
\end{equation} 
 
 If hypothesis E(R,f) is verified, there exists a subset of $A$ of large measure on which the above control holds pointwisely. More precisely,
  consider an arbitrary $\a\in (0,1)$ and  define the following subset of $A$:
  
 \begin{equation}\label{defB}
  B_\a:= \left \{x\in A:\quad {1\over 2}\Big[ \mathbf{P}(\tau^{ x}_{G}\le 2R)+\mathbf{P}(\tau^{\leftarrow x}_{G}\le 2R)\Big]\le f^{\a} \right \}. 
 \end{equation} 
 
The set $B_{\a}$ is almost of full measure within $A$, representing in some sense  the basin of attraction of 
the local equilibrium. Indeed using $E(R,f)$, one easily gets that
\be{piB}
\pi_A(B_\a)\ge 1-(f)^{1-\a}
\ee
since
$$f\ge f_A(2R)\ge {1\over 2} \sum_{x\in A\backslash B_\a }\pi_A(x)\Big[ \mathbf{P}(\tau^{ x}_{G}\le 2R)+\mathbf{P}(\tau^{\leftarrow x}_{G}\le 2R)\Big] \ge (f)^{\a}
\pi_A(A\backslash B_\a).$$ 
 \item[\textbf{T(R,d)}] \textbf{Fast thermalization in $B_\a$}: for $d \in (0,1)$, we say that $A$ satisfies the thermalization condition 
  $T(R,d)$ if  the process started from any point in $B_{\a}$ loses memory of its initial condition within a time scale $R$; more precisely
   \be{defT}
    \bar d_{B_{\a}}(R) \leq d.
    \ee     
\item[\textbf{Rc(R,r)}] \textbf{Fast recurrence to the set $B_{\a}\cup G$}: starting from anywhere inside of the state space, $X$
  reaches either $G$ or $B_{\a}$ with high probability within a time $R$; namely, let $r \in (0,1)$, we say 
  that $A$ satisfies the hypothesis $Rc(R,r)$ if 
   \begin{equation}\label{ricB}
\sup_{x\in \cX}\mathbf{P}(\t^{x}_{B_\a\cup G}>{R})<r.
     \end{equation}
        \end{enumerate}
        
         Note that both conditions $T(R,d)$ and $Rc(R,r)$ depend on the parameters $f$ and $\a$ since they depend on the set $B_{\a}$. Both for lightness 
         of notation and in view of Proposition \ref{mixgen}, which states that a control on the mixing time of the dynamics over the whole state space 
         coupled to $E(R,f)$ is sufficient to ensure that conditions $T(R,d)$ and $Rc(R,r)$ hold with explicit parameters, we
         choose not to write this dependence when refeering 
          to fast thermalization and to fast recurrence. 
       
        \begin{hyp}[Hypothesis HpG]
          When the set $A$ satisfies simultaneously the conditions $E(R,f)$, $Rc(R,r)$ and 
          $T(R,d)$, we say that condition $Hp.G(R,d,f,r)$ is verified. 
        \end{hyp}

In words, condition $Hp.G(R,d,f,r)$ for $d,f,r$ small means that the time scale $R$ is  large enough so that the process
loses memory of its initial condition in $B_{\a}$, but that within a time $2R$ the trajectories starting from 
the typical configurations for  $\p_{A}$ inside of $A$ still did not hit $G$, while hitting either $B_{\a}$
or $G$ within the same time scale occurs with high probability.  This is a typical situation characterizing metastability.
 
  Before stating our main result, a couple of remarks are in order: 

  \begin{enumerate}
  \item[$\bullet$]
The recurrence hypothesis $Rc$ can be verified by changing the time scale $R$.
Indeed, if $\tau_{G}^{\p_{A}}$ is much larger than  
$R$, our process satisfies hypotheses $E(R,f)$ and $T(R,d)$  within a time $R_+ \gg R$. In this case 
$Rc(R_+,r)$ holds with $r=\Big(\sup_{x\in \cX}\mathbf{P}(\t^{x}_{B_\a\cup G}>{R})\Big)^{\frac{R_+}{R}}$.
 \item[$\bullet$] 
  In view of equation \eqref{piB}, it would seem natural to replace the thermalization condition
  $T(R,d)$ by $\bar d_{A}(R) \leq d$. We stress that from a physical point of view, this is
  requiring too much on the whole set of trajectories starting from $A$ and that in most applications 
  there is a region (which is here the set $A \setminus B_{\a}$) in which some trajectories might escape 
  from $A$ within a short time without thermalizing afterwards. The typical picture of a well in
  which trajectories thermalize is here played by the basin of attraction $B_{\a}$.
  In different words, with these hypotheses we can apply our results also to cases in which
  the basin of attraction of the stable state is not known in details.
  
    Nevertheless, in the (degenerate) case where thermalization occurs over the whole set $A$,
    condition $Hp.G(R,d,f)$  can be verified. In this spirit, we state the following result. 
  \end{enumerate}
  
   \begin{proposition}\label{mixgen} 
    If $A$ satisfies hypothesis $E(R,f)$ and 
      \be{defTA}
   \bar d_{A}(R) \leq d,
    \ee     
    then $A$ satisfies condition $Hp.G\left(R,d,f,d+f+f^{1-\a} \right)$.  
    \end{proposition}

  We now state our main technical result.

\subsection{Convergence of  measures}\label{Luno}

Define   the following {\it doubly conditioned
 evolved measure on $A$}:
for any $t\ge 2R$ and $x\in A$ the evolution measure conditioned to $A$ in $[0,t]$   and to visit $B_\a$ within $t-2R$ is defined as:

\be{misdcond}
\hat \mu^x_{A,t}(.):=\mathbf{P}(X^x_t=. |\t^x_{G}>t, \t^x_{G\cup B_\a}\le t-2R)\ee
and
\be{conddpa}
\hat \mu^{\p_A}_{A,t}(.):=\sum_{x\in A} \pi_A(x)\hat\mu^x_{A,t}(.).
\ee

 \begin{remark}\label{hatB}
   For any $x\in B_\a$ we have $\hat\mu^x_{A,t}=\tilde \mu^x_{A,t}$ since $ \t^x_{G\cup B_\a}=0$
and:
\be{hatpa}
d_{TV}(\hat \m^{\p_A}_{A,t},\tilde\m^{\p_A}_{A,t})\le \p_A(A\backslash B_\a)<f^\a.
\ee
 \end{remark}

For $x\in A\backslash B_\a$ the measures  $\hat\mu^x_{A,t}$ and $\tilde \mu^x_{A,t}$ can be different,
even if the second conditioning $\t^x_{G\cup B_\a}\le t-2R$ is an event of large probability when $t>3R$.
In particular they can be different for $x$ and $t$ such that $\mathbf{P}(\t^x_{G}>t)<\mathbf{P}(\t^x_{G\cup B_\a}\le t-2R)$.

\bt{convdtv}
Under hypothesis  $Hp.G(R,d,f,r)$,  for every  $\a\in (0,1)$, for $d$ and $f$ such that $r+2f^\a<{1\over 4}$,
we define $\bar c:={1\over 2}-\sqrt{{1\over 4}-c}$.
 The following vicinity of the doubly conditioned evolution measure and the conditioned invariant one holds:  
\be{dtv}
\sup_{x\in A} \sup_{t\ge 2 R}d_{TV} \left(\hat\mu^{x}_{A,t},\pi_{A} \right)<4 \left[\bar c + 2f + f^{\a} + d\right]=:\e_1,
\ee
\be{dtvpa}
 \sup_{t>2R}d_{TV} \left(\tilde\mu^{\p_A}_{A,t},\pi_{A} \right)<4 \left[\bar c + 2f + f^{\a} + d\right] + f^{1-\a}=:\e_2=\e_1+f^{1-\a}.
\ee
 
 Moreover
\be{dtvqs}
\sup_{x\in A} \sup_{t\ge 2 R}d_{TV} \left(\hat\mu^{x}_{A,t},\m^*_{A} \right)<\e_1+\e_2.
\ee
\et
 
  Roughly speaking, Theorem \ref{convdtv} states that after time $2R$, the process thermalized in the metastable states in the sense of the doubly
   conditioned measure; for $x \in B_{\a}$, this result justifies
  the definition of metastable state for the conditioned evolved measure $\tilde \mu^{x}_{A,2R}$.  
   The important consequence of these results is the fact that the quasi stationary measures and the conditioned invariant ones are 
  asymptotically close in total variation.

The following is an immediate consequence of Theorem \ref{convdtv}:
\bc{cor1}
Under the same hypotheses then the ones of  Theorem \ref{convdtv}, for any $x\in B_\a$ we have
$$
d_{TV}\left(\m^{(em)x}_A,\p_A\right)<\e_1.
$$
\ec
Indeed for $x\in B_\a$ we have
$$
 \m^{(em) x}_A(y)=\frac{\mathbf{E}\xi^x_G(y)}{\mathbf{E}\t^x_G}=\frac{\int_0^\infty \tilde \m^x_{A,t}(y) \mathbf{P}(\t^x_G>t)dt}{\mathbf{E}\t^x_G}=
 \p_A(y)+\frac{\int_0^\infty (\tilde \m^x_{A,t}(y)-\p_A(y)) \mathbf{P}(\t^x_G>t)dt}{\mathbf{E}\t^x_G}.
$$

\subsection{
Exponential behaviour}

 A consequence of the results of the previous section is the 
following asymptotic exponential behavior:
 
 \begin{theorem}\label{expbeh} 
  Assume that hypothesis $HpG$ is satisfied. There exists $\e$ such that $\e = O(r+ \e_{2})$ where $\e_{2}$ has been
  defined in Theorem \ref{convdtv}
   such that 
  the following inequality holds for any $t \geq 0$:
    \begin{equation}\label{fortispiAl}
        \left| \mathbf{P}\left[ \frac{\t^{\pi_{A}}_{G}}{\mathbf{E}[\t^{\m^*_{A}}_{G}] } \geq t \right] - e^{-t} \right| \leq  \e e^{-t}.  
       \end{equation} 
       
        A pointwise
        control also holds starting from $B_{\a}$; namely, for every $t \geq 0$:
    \begin{equation}\label{fortisl}
      \sup_{x \in B_{\alpha}} \left| \mathbf{P}\left[ \frac{\t^{x}_{G}}{\mathbf{E}[\t^{\m^*_{A}}_{G}] } \geq t \right] - e^{-t} \right| \leq  \e e^{-t} 
    \end{equation} 
    
    Moreover, for $x\not\in B_\a$:
     \begin{equation}\label{fortisB}
       \left| \mathbf{P}\left[ \frac{\t^{x}_{G}}{\mathbf{E}[\t^{\m^*_{A}}_{G}] } \geq t \right] - \mathbf{P}(\t^x_G>2R,
    \t^x_{G\cup B_\a}\le R)
    e^{-t} \right| \leq  \e e^{-t}.
     \end{equation} 
  
  \end{theorem}
  
   Let us mention some straightforward consequences of the strong control on hitting times provided by Theorem \ref{expbeh}. The first one 
    states an equivalence of time scales in our setup, and directly follows from integrating the relations \eqref{fortispiAl} and 
   \eqref{fortisl} for $t \in \mathbf{R}^{+}$. 
    
      \begin{corollary}
        For any $x \in B_{\a}$, one has 
       \begin{equation}\label{eqT}
        \mathbf{E}[\t^{\m^*_{A}}_{G}] \sim  \mathbf{E}[\t^{\pi_{A}}_{G}] \sim \mathbf{E}[\t^{x}_{G}].
       \end{equation} 
      
      \end{corollary}

     The second one 
    is an extension of \eqref{fortispiAl} to a wide set of starting measures.  
   
    \begin{corollary} For any probability measure $\mu$ such that $Supp(\mu) \subset B_{\a}$ and any $t \geq 0$
         \begin{equation}\label{fortispiAmu}
        \left| \mathbf{P}\left[ \frac{\t^{\mu}_{G}}{\mathbf{E}[\t^{\mu}_{G}] } \geq t \right] - e^{-t} \right| \leq  \e e^{-t}.  
       \end{equation}
        
         In particular, the equivalence \eqref{eqT} also generalizes to  any probability measure $\mu$ such that $Supp(\mu) \subset B_{\a}$, namely
        \begin{equation}
        \mathbf{E}[\t^{\m^*_{A}}_{G}] \sim  \mathbf{E}[\t^{\mu}_{G}].
       \end{equation} 
    \end{corollary}

  \begin{remark}\label{remdop}  The following  observations are in order:
  \bi
  \item[-] The setting of Theorem \ref{expbeh} is quite general; in particular, it is a non reversible setup, and it is noteworthy
   to realize that the invariant measure of the set of goals $\p(G)$ can be large. Roughly speaking, the goal set does not need 
    to be a \textit{rare event}.
     
  \item[-] The exponential control on the queues of distribution is exactly of order $e^{-t}$ for $t \geq 0$; in particular, 
   we do not lose on the exponential controlling the error as in \cite{OI}, where the control on the queue is of the type 
    $e^{-(1-c)t}$ where $c$ is a small constant not depending on $t$.
     
  \item[-] Our results hold pointwise for a large set within the complement of $G$, and yield an explicit control on the errors 
   with respect to the hitting time. 
    
  \item[-] For  $x\in A \setminus B_\a$, the picture is intuitively the following:  in a time of order $R$, the process
  either thermalizes in $B_\a$ or hits the set $G$; in the first case, by Markov's property, we are back to the case of a process with starting point 
   within $B_{\a}$. 
  \ei 
  
  \end{remark}

  For practical purposes, we stress that a consequence of Theorem \ref{expbeh} can be deduced in the asymptotic regime of Section \ref{Not}; 
  more precisely, 
   assume that hypotheses $E(R_{n},f_{n})$, $Rc(R_{n},r_{n})$ and $T(R_{n},d_{n})$ hold for a family of Markov processes 
   $(\cX^{n},X^{n})$ with 
    parameters $f_{n}, r_{n}$ and $d_{n}$ which are $o(1)$; then the following result holds: 
  
   \begin{corollary}\label{convlaw}
       Under hypothesis $Hp.G_{n}$, as $n \to \infty$, the following convergence in law holds 
   \begin{equation}\label{glconvl}
     \t^{\pi_{A}}_{G^{(n)}}/\mathbf{E}[\tau^{\pi_{A}}_{G^{(n)}}]  \stackrel{\cL}{\rightarrow} \mathcal{E}(1)
    \end{equation} 
         where $\mathcal{E}(1)$ denotes a random variable which is exponentially distributed with mean one.
          
    Furthermore, for any starting measure $\mu^{(n)}$ with support contained in $ B_{\alpha}^{(n)}$, the following convergence in law holds:  
    \begin{equation}\label{ptconvl}
     \t^{\mu^{(n)}}_{G^{(n)}}/\mathbf{E}[\tau^{\mu^{(n)}}_{G^{(n)}}]  \stackrel{\cL}{\rightarrow} \mathcal{E}(1).
    \end{equation} 
   \end{corollary}

\section{Applications }\label{Ex}

In this part we discuss two examples. 

 Our first example originates from a model for dynamical polymers  which was first introduced and analyzed in \cite{CMT}. The asymptotic
  exponentiality of the tunneling time was shown in \cite{CLMST} by considering in particular the death and birth process of our first example.

The second case is a model of shuffling cards,  well known as "top-in-at-random" (TIAR).
This is actually a case where the invariant measure is the uniform one and
the "potential barrier" is only entropic. 

\subsection{Birth and death processes with logarithmic barrier}
\label{sec3.1}

 {\bf The model}
 
  In this part, $c > 0$ is a constant which may vary from line to line.
 
  The (typically unique) crossing point of the interface of the dynamical polymer model with the wall in the delocalized phase follows the 
   law of the  birth and death process on $\{0,...,n\}$
 with invariant measure
 $$
 \pi(x):=Z^{-1}{1\over (x\vee 1)^{3/2}((n-x)\vee 1)^{3/2}}
 $$
 where $Z=Z(n)$ is a normalization constant. We refer to \cite{CLMST}[Section 1.2] for details.
 
  We note that it is easy to show that
   \begin{equation}\label{eqZ}
    Z(n) \sim c n^{3/2}.
   \end{equation} 
  
 We consider the Metropolis dynamics associated to this invariant measure, which means that the birth rate $b(\cdot)$ of the chain is given by
  \begin{equation}\label{bir}
   b(x)=\min \left \{1,{\pi(x+1)\over \pi(x)} \right \},
   \end{equation}  so that 
 $Q(x,x+1)=b(x)$. Moreover, the death rate $d(\cdot)$ is given by
 \begin{equation}\label{dea}
  Q(x,x-1)=d(x)=\min \left \{1,{\pi(x-1)\over \pi(x)} \right \}.
  \end{equation}

  This birth and death process has been introduced in \cite{CLMST} to describe the evolution of the (typically unique) crossing
  location of the pinning line by a polymer chain in the delocalized regime. 
 
From \eqref{bir} and \eqref{dea}, one readily observes that for $x\le n/2$, the process exhibits a small drift towards $\{0\}$ since $b(x)<1$
 and $d(x)=1$,
while for $x>n/2$ there is a small drift towards $\{n\}$ since $b(x)=1$ and $d(x)<1$. Otherwise stated, we
are considering a one dimensional random walk under a potential given by
 a symmetric double well with minima in $0$ and $n$ and with a small logarithmic barrier at $\{n/2\}$.

It is shown in \cite{CLMST} that the mean time for the process starting from $\{0\}$ to reach $\{n/2\}$, and hence $\{n\}$,  is 
of order $n^{5/2}$;
however, one directly observes that the ratio ${\pi(0)\over \pi(n/2)} = \frac{n^{3/2}}{8}$. This entails that this model
is a case where the metastable behavior
is much more delicate with respect to the Metropolis case where the tunneling time
is of order of the ratio of the value of the invariant measure between  the bottom and
 the top of the barrier, see \cite{OV}[Chapter 6,Theorem 6.23] for details about this well known phenomenon; as a consequence,
 rough arguments based on reversibility have to fail in this model, see
in particular point ii) in \cite{OV}[Chapter 6,Theorem 6.23].
 The main result of this part is the following: 
 
 \begin{theorem}\label{BaD} Define $G := \{ n/2,n/2+1,\ldots, n\}$, $\a \in (0,1)$ and $n$
 large enough such that $n^{\a} < n/2$. HpG holds with $n^{5/2-\e} << R  $ for small enough $\e$.

  \end{theorem} 

 \textit{Proof}

   Our aim is to exhibit a sequence $(R_{n})_{n \geq 0}$ such that the hypothesis of Theorem \ref{expbeh} are satisfied. 
   In fact, we will show that the stronger hypothesis of Proposition \ref{mixgen} actually hold.
     We define the set 
    \begin{equation}
    B_{\a} := \{ 0, \ldots, n^{\a} \}. 
    \end{equation}
   
    Our first step is to show that for any $x\in\{0,...,n/2-1\}$, the following equivalence holds as $n \to \infty$:   
     \begin{equation}\label{tp}
     \mathbf{P}(\t^{x}_0>\t^{x}_{n/2})\sim \left ({2x\over n} \right )^{5/2}.
     \end{equation}
     
      We first recall the classical notion of resistance between two states $x$ and $y$:
      \begin{equation}
      R(x,y) = \left ( \pi(x) P(x,y)\right )^{-1}.
      \end{equation}
       
   In this one--dimensional setup, it is a classical result that the resistances are linked to the law of first hitting to a set through 
    the basic relation (see for example \cite{PeLi}[Chapter 2])  
   \begin{equation}
 \mathbf{P}[ \tau^{x}_{0} > \tau^{x}_{y}] = R(x,y)
   \end{equation} 
   where $x < y$ and $R(x,y)$ is obtained by summation from the $R(k,k+1)$, namely 
   \begin{equation}
    R(x,y) = \sum_{k=x}^{y-1} R(k,k+1).
   \end{equation} 
    
    In our specific case, it is an easy computation to realize that, as $x \to \infty$:
    \begin{equation}\label{eqR}
     R(x,x+1)=\Big({\p(x)\over 1+({x\over x+1})^{3/2}}\Big)^{-1}\sim x^{3/2}+(x+1)^{3/2}.
     \end{equation} 
      
Indeed, representing the Markov chain like a series of resistances, we deduce that
$$
R(0,x)\sim cx^{5/2}
$$
and the probability in (\ref{tp}) representing the potential at point $x$ is given by
$$
 \mathbf{P}(\t^{x}_0>\t^{x}_{n/2})= V(x)={R(0,x)\over R(0,n/2)}\sim \left({2x\over n}\right)^{5/2}.
$$
 
From \eqref{tp}, one deduces that there exists a constant $c > 0$ such that
 \begin{equation}\label{tp1}
 \mathbf{P}(\t^{n^\a}_0<\t^{n^\a}_{n/2})\ge 1- c\big(2n^{\a-1}\big)^{5/2}.
 \end{equation}
  
 Then we define two sequences of stopping times $(\tau_{i})$ and $(\sigma_{i})$ in the following way:  
$$
\t_0:=\inf\{t:\; X_t=0\},\quad \s_i:=\inf\{t>\t_{i-1}:\; X_t=n^\a\}, \quad \t_i:=\inf\{t>\s_{i}:\; X_t=0\}
$$
and, for a given $T > 0$, we introduce 
$$
\nu(T):= \max \left \{i:\; \t_i<T\wedge \t^{n^\a}_{n/2} \right \}.
$$
 
For any $x\le n^\a$ and $N>1$, making use of \eqref{tp1}, we can write 
 \begin{equation}\label{hitn2} \begin{split}
  \mathbf{P}(\tau^{ x}_{n/2}> T)\ge \mathbf{P}(\tau^{ x}_{n/2}> T,\; \nu(T)<N) & =
\mathbf{P}(\tau^{ x}_{n/2}> T\;|\; \nu(T)<N)\mathbf{P}( \nu(T)<N) \\
 & \ge \Big(1- \big(2n^{\a-1}\big)^{5/2}  \Big)^N
\mathbf{P}( \nu(T)<N).
  \end{split}
 \end{equation}
  
On the other hand, we can show that 
 \begin{equation}\label{Enu}
 \mathbf{E}(\nu(T))\le {T\over \mathbf{E}(\t^{0}_{n^\a})}.
 \end{equation}

Equation \eqref{Enu} is obtained in the following way: we consider the martingale $(M_{k})$ defined by
 \begin{equation}
 M_k :=\sum_{i=1}^k 
\Big[\t^{0}_{n^\a}- \mathbf{E}(\t^{0}_{n^\a})\Big]. 
 \end{equation}
  
 Since 
  $T\ge {\nu}(T) \t^{0}_{n^\a}$, we get the inequality
  \begin{equation}
  M^{\nu(T) }\le T-
\nu(T) \mathbf{E}(\t^{0}_{n^\a}).
  \end{equation}
   
  Then we apply Doob's optional-stopping theorem to $M$ at time $\nu(T)$, and \eqref{Enu} follows. 
  
   Now we prove the equivalence:
  \begin{equation}\label{hit0}
   \mathbf{E}(\t^{0}_{n^\a}) \sim c n^{ 5\a/2}.
   \end{equation} 

   We make use of the equivalences \eqref{eqR} and \eqref{eqZ} and we write:
    \begin{equation}\label{pot}
     \begin{aligned}
       & \mathbf{E}(\t^{0}_{n^\a}) = \sum_{k=0}^{n^{\a}-1} \pi(k) R_{n^{\a}}^{k} =  \sum_{k=0}^{n^{\a}-1} \pi(k) \sum_{i=k}^{n^{\a}-1} R_{i}^{i+1} \\
    & \phantom{iii} \sim  c \sum_{k=0}^{n^{\a}-1} \pi(k) \sum_{i=k}^{n^{\a}-1} i^{3/2} \sim c \sum_{i=0}^{n^{\a}-1} \sum_{k=0}^{i}
      i^{3/2} \frac{n^{3/2}}{ (1 \vee k)^{3/2} (n-k)^{3/2} }.
     \end{aligned}
    \end{equation}

 Since uniformly on $k \in [0,n^{\a}]$, we have
\begin{equation}
 \frac{n^{3/2}}{(n-k)^{3/2}} \sim 1,
\end{equation} 
  we deduce from \eqref{pot} that 
  \begin{align}
   \mathbf{E}(\t^{0}_{n^\a}) \sim c \sum_{i=0}^{n^{\a}-1} i^{3/2} \sum_{k=0}^{i} \frac{1}{(1 \vee k)^{3/2}} \sim  c \sum_{i=0}^{n^{\a}-1} i^{3/2} 
   \sim c n^{5\a/2},
  \end{align}
 and hence we recover \eqref{hit0}.
  
  Now we prove that the conditions of Proposition \ref{mixgen} hold. 
  
  \textit{ Mixing condition}

Fix $\e>0$ sufficiently small and consider $R=R_n=n^{{5\over 2}-\e}$ and $N=n^{(1-\a){5\over 2}-{\e\over 2}}$.
We first  prove that $d_{A}(R_n)\to 0$ by using its estimate in terms of coupling time (
see for example \cite{LPW}[Chapter 5, Theorem 5.2] for this classical result):
 \begin{equation}\label{spieg}
 d(R_n)\le  \mathbf{P}(\t^{n/2}_0>R_n)\le {\mathbf{E} \t^{n/2}_0\over R_n}.
 \end{equation}
  
 We show the bound
  \begin{equation}\label{rgh}
  \mathbf{E}\left[ \t^{n/2}_0\right] \leq c n^{2}.
  \end{equation}
   
  Indeed, we have the inequalities:
  \begin{align*}
   \mathbf{E}[ \t^{n/2}_0] = \sum_{k=0}^{n/2} \pi(k) R_{0}^{k} \leq c \sum_{k=0}^{n/2} \frac{1}{(1\vee k)^{3/2}} \sum_{i=0}^{k} i^{3/2} \\
    \leq c \sum_{i=0}^{n/2} i^{3/2}  \sum_{k=i}^{n/2} \frac{1}{(1\vee k)^{3/2}} \sim c \sum_{i=0}^{n/2} i^{3/2} i^{-1/2} \sim c n^{2},
  \end{align*} 
  from which \eqref{rgh} follows. 
   
 Hence the quantity in the right hand side of \eqref{spieg} vanishes as $n \to \infty$, and \eqref{defTA} holds 
  with $R = R_{n}$. 
  
   \textit{Slow escape}
  
On the other hand, combining Markov's inequality and \eqref{hit0}, we get: 
 \begin{equation}\label{hitso}
 \mathbf{P}( \nu(T)<N)\ge 1-{\mathbf{E}(\nu(T))\over N}\ge 1- c{T\over n^{ 5\a/2} N}=1-cn^{-\e/2}. 
 \end{equation}

 Finally we combine \eqref{hitn2} and \eqref{hitso} to deduce that, for any $x \in \cB_{\a}$:
 \begin{equation}
 \mathbf{P}(\tau^{ x}_{n/2}> T)\ge(1-c n^{-\e/2})(1-n^{-\e/2})\sim 1-Cn^{-\e/2}.
 \end{equation}
We finally conclude
$$
\sum_x\p(x)\mathbf{P}(\tau^{(n), x}_{n/2}\le 2 T_n)\le \sum_{x\le n^\a}\p(x)Cn^{-\e/2}+\sum_{x> n^\a}\p(x)
\le Cn^{-\e/2}+C'n^{-\a/2}\to 0,
$$
which achieves the proof of $E(R_{n},f)$; thus we can apply Theorem \ref{expbeh} to deduce Theorem \ref{BaD}. 

 $\qed$

\subsection{Top-in-at-random, an entropic barrier}\label{sec3:2}

   In this part, we apply the results of the section \ref{Luno} to the Top-in-at-random model. 
   The state space is $\cX^{(n)} = \{ \text{permutations of } (1,\ldots,n) \} $
   and it is called the {``deck of $n$ cards"}. The dynamics is the following: at each step, we take the upper card and put it in a random position; more precisely, we update the configuration 
    $x = (x_{1}, \ldots, x_{n})$ to a new  configuration
    $x^k = ( x_{2}, \ldots, x_{k}, x_{1}, x_{k+1},\ldots, x_{n})$, with $k$ uniform in $\{1,\ldots,n\}$. 
    We call the transition form $x$ to $x^k$ a ``shuffle". 
    By symmetry, it is straightforward to  realize that the invariant measure of this dynamics is uniform  on $\cX^{(n)}$. 
     
      This model is well known, it is historically the first model which was shown to exhibit cutoff (\cite{AlD}). 
      Here we will only use the fact that in the same paper the authors show that $R_{n} \sim n \log(n)$.

    We consider $G := \{ (1,\ldots,n )\}$ the ordered configuration, and we define a projection 
     \begin{equation}
      \s(x) = \max \{ i \leq n, \; x_i > x_{i+1} \}.
     \end{equation} 
    
    It is easy to see that $\s(x)$ measures the minimal number of 
    shuffles 
      to reach $G$ starting from $x$; in particular, $\s(x) = 0$ if and only if $x = G$. 
      
\bl{TIARproj}      The process $\s(X^{x}_{t})$
       is still Markovian with transition probabilities given by:
\begin{equation}\label{proj}
       P_{i,j} :=
       \mathbf{P}\left[ \s(X_{t+1}^{x}) = j | \s(X_{t}^{x}) = i  \right] =   \left \{  
       \begin{array}{lll}
       1/n & \text{if} & j = i-1 \\
       i/n & \text{if} & j = i \\
       1/n & \text{if} & j = i+1, \dots,n\\
       \end{array} \right .
      \end{equation} 
\el
\emph{Proof:} 

Let $y=X^x_t$ be the current configuration and 
$k$ be the new position for the first card, so that $X^x_{t+1}=y^k$. 
Let $\s:= \s(y)$.
If the first card is inserted above $\s$, namely if $k < \s$, then $\s(y^k)= \s$.
There are $\s -1$ available choices.

Otherwise, consider the stack formed by the last $n-\s$ cards of $y$ plus the card $y(1)$.
Each of the possible positions $k\in \{\s ,\ldots,n \}$ corresponds to a different 
$\s (x^k)\in \{\s-1,\ldots,n-1 \}$. 
To see this, it is sufficient to observe that 
$k\longrightarrow \s (x^k)$ is invertible in  
$\{\s ,\ldots,n\}$, the inverse application being $\{\s-1,\ldots,n-1 \}$.
Let  
$k^*:=\max\{i \in \{\s , \ldots,n \} \ ; \; y(i) < y(1)\}$ 
be the ``natural position" of the card $x(1)$ into this stack (notice that $k^*=\s$ if $y(1)<y(\s +1)$ ). 

We get the value 
$\s (y^k)=\s'\in \{\s -1,\ldots,n-1\}$ 
by choosing
\be{kdisigma}
  k(\sigma')=
   \begin{cases}
    k^{*} & \text{ if }\sigma'=\sigma-1,\\
    \sigma' & \text{ if }\sigma-1<\sigma'<k^{*},\\
    \sigma'+1 & \text{ if }\sigma'\ge k^{*}.
   \end{cases}
\ee
 
Thus, under the condition $k \ge \s $,
$ \s(X^x_{t+1})$ takes the values in $\{ \s-1,\ldots,n-1\}$
with uniform probability. 

$\qed$

We denote by $\s^i(t)$ the Markov process with transition matrix $P$ defined in (\ref{proj}).
      


 To show that $HpG$ holds, we need to get good estimates on the hitting times $\tau^{k}_0$ of $\{0\}$ for the projected chain $\s^{k}$. 
 
 
 Let $\xi^{i}_{j}(k) := \left| \{ t \leq \tau^{i}_{j}, \s^i(t)=k   \}\right|$ be the local time spent in $k$ before
 hitting $j$ when starting from $i$.
 We define $\xi^i_j:=\xi^i_j(i)$.
 
 Since in the downward direction only one-step transitions are allowed, we have: 
  \begin{equation}\label{beg}
  \begin{aligned}
    \mathbf{E}\left[ \tau_{0}^{n-1} \right] &=
  \mathbf{E}\left[ \sum_{k=1}^{n-1} \xi_{0}^{n-1}(k) \right]\\ 
  &=\sum_{k=1}^{n} \mathbf{E}\left[ \xi_{0}^{k} \right]
  \mathbf{P}\left[ \t^{n-1}_{k}<\tau^{n-1}_0 \right] \\
  &=
  \sum_{k=1}^{n} \mathbf{E}\left[ \xi_{0}^{k} \right],
  \end{aligned}
  \end{equation} 
  where we used the strong Markov property at time $\tau^0_k$ in the second identity. 
    
Since
 \begin{equation}\label{expxi}
  \begin{aligned}
    \mathbf{P}[ \xi_{j}^{i} > n] 
    &=
    \mathbf{P}[ \tau_{i}^{i} < \tau_{j}^{i}] 
    \mathbf{P}[ \xi_{j}^{i} > n-1]
    \\
    &=
    \mathbf{P}[ \tau_{i}^{i} < \tau_{j}^{i}]^n,
  \end{aligned}
 \end{equation} 
$\xi^i_j$ is a 
geometric variable with mean 
$\mathbf{P} [\tau_{j}^{i} < \tau_{i}^{i}]^{-1}$.

Let $I$ be a subset of $\{0,\ldots,n-1\}$.   We recall the standard ``renewal" identity 
  \begin{equation}\label{renew}
   \mathbf{P}[ \tau^{i}_{I} < \tau_{j}^{i}] = \frac
   {\mathbf{P}[ \tau^{i}_{I} < \tau_{j \cup i}^{i}]}
   {\mathbf{P}[ \tau_{I\cup j}^{i} < \tau_{i}^{i}]}. 
  \end{equation} 

Our main estimate is contained in the following Lemma. 
\bl{xik0}  For any $k \in [1,n]$,
\be{stimaxik0} 
  \mathbf{P}[\tau^k_0<\tau^k_k]
  \le
  \frac{(n-1)(n-k-1)!}{n!}.
\ee
\el
\emph{Proof:}

By \eqref{renew},
 \begin{equation}\label{stimaxi}
  \begin{aligned}
   \mathbf{P}\left[ \tau_{0}^{k} < \tau_{k}^{k} \right]  
      & = P_{k,k-1} \mathbf{P}\left[  \tau_{0}^{k-1} < \tau_{k}^{k-1} \right] \\ 
     & = P_{k,k-1} 
     \frac
     {
     \mathbf{P}\left[ \tau_{0}^{k-1} < \tau_{k,k-1}^{k-1} \right]}
     {
     \mathbf{P}\left[ \tau_{0,k}^{k-1} < \tau_{k-1}^{k-1} \right]}. 
  \end{aligned}
 \end{equation}

  Since the downward moves are single-step, for $k \geq 1$, we have:
   \begin{equation}\label{single}
    \mathbf{P}\left[ \tau_{0}^{k-1} < \tau_{k,k-1}^{k-1} \right] = \mathbf{P}\left[ \tau_{0}^{k-1} < \tau_{k-1}^{k-1} \right]. 
   \end{equation}  
    
On the other hand, 
 \begin{equation}\label{denomxi}
  \begin{aligned}
   \mathbf{P}\left[ \tau^{k-1}_{0,k} < \tau_{k-1}^{k-1}  \right] & =  
     \mathbf{P}\left[  \tau^{k-1}_{0} < \tau_{k-1}^{k-1}  \right] + \mathbf{P}[ \tau^{k-1}_{k} < \tau_{k-1}^{k-1} < \tau^{k-1}_{0}  ]\\
        & \ge   \mathbf{P}\left[  \tau^{k-1}_{0} < \tau_{k-1}^{k-1} \right] + 1 - P_{k-1,k-2} - P_{k-1,k-1}
  \end{aligned}
 \end{equation} 
where we used again the fact that downhill moves are unitary.

Let 
$E_k:=\mathbf{P}\left[ \tau^{k}_{0} < \tau_{k}^{k}  \right]^{-1}$. 
Plugging (\ref{single}), (\ref{denomxi}) into (\ref{stimaxi}), and using (\ref{proj}) we get

  \begin{equation}\label{recxi}
   E_k \ge n + (n-k)E_{k-1}.  
  \end{equation} 
   
    Since $E_1 = p_{1,0}^{-1}=n$, we prove inductively  that    
     \begin{equation}\label{solrec}
      E_k \ge n \sum_{j=1}^k\frac{(n-j-1)!}{(n-k-1)!}.  
     \end{equation} 
      
Indeed, (\ref{solrec}) holds for $k=1$ and, by (\ref{recxi}), (\ref{solrec}),
\begin{equation}
 \begin{aligned}
   E_{k+1} 
  &\ge 
  n\left(1 + (n-k-1)\sum_{j=1}^k\frac{(n-j-1)!}{(n-k-1)!}\right) \\
  &=
  \frac {n}{(n-k-2)!} \left((n-k-2)! + \sum_{j=1}^k {(n-j-1)!}\right).
 \end{aligned}
\end{equation}

By taking only the largest term in (\ref{solrec}), \eqref{stimaxik0} immediately follows.

$\qed$

\bl{TIARmu}
The invariant measure of $\s^{i}$ is 
\be{muTIAR}
 \bar \mu(k)
 =
 \begin{cases}
 \frac{1}{n!} & \text{ if } \   k=0, \\
 \frac{n-k}{(n-k+1)!} & \text{ if }  \ k>0.
 \end{cases}
\ee
\el
\textit{Proof} 

%
%
Direct computation shows that,
for $k\in\{0,\ldots,n \}$:
\be{muTIAR1}
  \bar \mu(\{0,\ldots,k\})=\frac{1}{(n-k)!},
\ee
and therefore
 \begin{equation}\label{TIARinvmu}
  \begin{aligned}
    \sum_{i=0}^{n} \bar \mu(i) P_{i,j}
  &=
  \frac 1 n \sum_{i=0}^{j+1} 
  \bar \mu(i)  
  +
  \frac {j-1} n  \bar \mu (j)  \\
  &=
  \frac 1 n 
  \frac{1}{(n-j-1)!}
  +
  \frac {j-1} n  \bar \mu(j) \\   
  &=
  \frac 1 n
  \left( 
  \frac{1}{(n-j-1)!}   
  \frac{(n-j+1)!} {n-j}  
  +
  (j-1)   
  \right)
  \bar \mu(j) \\     
  &=
  \frac{(n-j+1)+(j-1)}{n}
  \bar \mu (j)
  =  
  \bar \mu (j). 
  \end{aligned}
 \end{equation} 

$\qed$   

\bt{TIARG+}
HpG holds in the TIAR model with $n \log n \ll R_{n} \ll (n-2)! $.
\et
\emph{Proof:}

 We get: 
 \begin{equation}\label{TiarG+}
   \begin{aligned}
    \sum_{x \in \cX_n} \mu(x)  \mathbf{P}(\tau^{ x}_{G} < R) 
  &=
  \sum_{k =0}^{n-1} \bar \mu(k)  \mathbf{P}(\tau^{ k}_{0} < R) \\
  &\le  
  \sum_{k =1}^{n-1} \bar \mu(k)  \mathbf{P}(\xi^{ k}_{0} < R) + \bar \mu(0) \\
  &=  
  \sum_{k =1}^{n-1} \bar \mu(k)  \left(1-\mathbf{P}[ \tau_{k}^{k} < \tau_{0}^{k}]^R\right)+ \bar \mu(0) \\
  &\le 
    \sum_{k =1}^{n-1} \bar \mu(k) R \; \mathbf{P}[ \tau_{0}^{k} < \tau_{k}^{k}]+ \bar \mu(0),
   \end{aligned}
 \end{equation}  
where we used (\ref{expxi}) in the third identity, combined with the fact that for $r >1$ and for any $x \in (0,1)$, one has the inequality
\begin{equation}
 1-x^{r} \leq r(1-x).
\end{equation} 

By using Lemma \ref{xik0} and Lemma \ref{TIARmu}, we get 
\begin{equation}\label{TiarG+'}
\begin{aligned}
 \text{r.h.s of (\ref{TiarG+})} 
  &\le
  R \sum_{k =1}^{n-1} 
  \frac{n-k}{(n-k+1)!} 
  \frac{(n-1)(n-k-1)!}{n!} +\frac 1 {n!}\\
  &= \frac {R(n-1)} {n!} \sum_{j=2}^n \frac 1 j + \frac 1 {n!} \\
  &= \frac {R \; n \log n} {n!} (1+o(1)).
\end{aligned}
\end{equation}

$\qed$

\section{ Proofs of Section \ref{Luno}}\label{Pro}

 \subsection{Control on evolution measures for $t\le 2R$}
 
  We first prove that, starting from $\p_{A}$, the law of $X$ stays close to  $\p_{A}$ at least until time $2R$: 
  
  \begin{lemma}\label{deleq} 
    If $E(R,f)$ holds, then, for $t \leq 2R$
     \begin{equation}\label{d12}
     d_{TV}(\m^{\p_A}_t,\p_A)\le f.
     \end{equation}
   \end{lemma}
   
    \textit{Proof} 
    
   We make use of the invariance of $\pi$ and we write
$$
d_{TV}(\m^{\p_A}_t,\p_A)={1\over 2\p(A)}\sum_{y\in\cX}\Big| \sum_{x\in A}\p(x)\mathbf{P}(X^x_t=y)-\p(y) \ind_{\{ y\in A\}} \Big|
$$
$$
={1\over 2\p(A)}\sum_{y\in\cX}\Big| \sum_{x\in A}\p(x)\mathbf{P}(X^x_t=y)- \sum_{z\in \cX}\p(z)\mathbf{P}(X^z_t=y) \ind_{\{ y\in A\}} \Big|
$$
$$
\leq {1\over 2\p(A)}\Big[\sum_{y\in G} \sum_{x\in A}\p(x)\mathbf{P}(X^x_t=y)+\sum_{y\in A} \sum_{z\in G}\p(z)\mathbf{P}(X^z_t=y)\Big].$$
 
Recalling the definition of the time reversal process $X^{\leftarrow}$
$$
\p(z)\mathbf{P}(X^z_t=y)=\p(y)\mathbf{P}(X^{\leftarrow y}_t=z),
$$
and noting that $\sum_{y\in G}\mathbf{P}(X^x_t=y)\le\mathbf{P}(\t^x_G\le t)$ (and that the same inequality holds for 
the time reversal process), we get
$$
d_{TV}(\m^{\p_A}_t,\p_A)\le f_A(t)\le f
$$
as soon as  $t\le 2R$. 

 $\qed$
  
\bl{mis}
If $E(R,f)$  and $T(R,d)$ hold, then for any $t\in[R,2R]$,
we have
\be{agg1}
\sup_{x\in B_\a}d_{TV}(\m^{x}_t,\m^{\p_A}_t)\le d+f^{1-\a}.
\ee
\el

 \textit{Proof} 
 
  We directly get: 
 \be{agg2}
\sup_{x\in B_\a}d_{TV}(\m^{x}_t,\m^{\p_A}_t)=\sup_{x\in B_\a}\frac{1}{2}
\sum_{y\in\cX}\Big| \sum_{z\in A}\frac{\p(z)}{\p(A)} \mathbf{P}(X^z_t=y)- \mathbf{P}(X^x_t=y) \Big|
\ee
$$\le \sup_{x\in B_\a}\sup_{z\in B_\a}\frac{1}{2}\sum_{y\in\cX}\Big|  \mathbf{P}(X^z_t=y)- \mathbf{P}(X^x_t=y) \Big|+\frac{\p(A\backslash B_\a)}{\p(A)}
\le d+f^{1-\a}.
$$
 $\qed$

\textit{Proof of Proposition \ref{mixgen}}   
 
  Let us assume now that $\bar d_{A}(R) \leq d$. We deduce 
\be{l13}
\sup_{x\in A} d_{TV}(\m^{\p_A}_t,\m^x_t)=
\sup_{x\in A}{1\over 2}\sum_{y\in\cX}\Big| \sum_{z\in A} \frac{\p(z)}{\p(A)}\mathbf{P}(X^z_t=y)- \mathbf{P}(X^x_t=y)\Big|\le\bar d_A(t)
\le d
\ee
for $t=R$. By triangular inequality, combining Lemma \ref{deleq} and equality \eqref{l13}, we obtain that 
 \begin{equation}
  \sup_{x\in A} d_{TV}(\p_A,\m^x_R) \leq d+f.
  \end{equation} 
  
 Hence we get that
  $$
\sup_{x\in\cX}\mathbf{P}(\t^{x}_{B_\a\cup G}>R)= \sup_{x\in A}\mathbf{P}(\t^{x}_{B_\a\cup G}>R)
$$
\begin{equation}\label{l1}
 \le \sup_{x\in A}\left(1-\mathbf{P}(\t^{x}_{B_\a}\le R) \right) \le \sup_{x\in A}(1-\m^x_R(B_\a))
\le d+f+f^{1-\a},
\end{equation} 
which proves Proposition \ref{mixgen}.  $\qed$

\subsection{Proof of Theorem \ref{convdtv} }

The main point is the proof of \eqref{dtv}. Indeed, \eqref{dtvpa} follows from noting that \eqref{dtv} implies
$$
\sup_{t>2R} d_{TV}(\hat\m^{\p_A}_{A,t},\p_A)<\e_1,
$$
 and then we combine \eqref{hatpa} and the triangular inequality.

 On the other hand, (\ref{dtvqs}) will follow from \eqref{dtvpa}, from the triangular inequality and from  the convergence 
$$
\m^*_A(y)=\lim_{t\to\infty} \tilde\m^{\p_A}_{A,t}(y),
$$
 which, from \eqref{dtvqs}, implies
 $$
d_{TV}(\p_A,\m^*_A)<\e_2.
$$

The strategy to prove \eqref{dtv} is to use the recurrence in the set $ B_\a$; let us first 
define for each $ s\ge 0$ and $x \in \cX$ 
 \be{st}
 \s^x(s):=\inf\{s'>s:\, X^{x}_{s'}\in B_\a\cup G\}.
 \ee

 We will follow the general guideline:
\bi
\item[1.] we first control the distance between the conditioned and the restricted measures 
  for starting points $z\in B_\a$ and $t\in[R,2R]$;
\item[2.] then we prove  estimates on the distribution of $\t^{x}_{G} $ in Lemma
\ref{l0};
\item[3.]  we  complete the proof combining the previous ingredients and strong Markov property at
time $\s(t-2R)$;
 \item[4.] finally we complete the proof for every starting point in $A$. 
\ei
Note that these ingredients  are frequently used in metastability in order to control the loss of memory with respect to
 initial conditions
 and consequently to deduce exponential behavior. 
We refer to \cite{FMNS} for similar results in finite state space regime, where the basin $B_{\a}$ is replaced
 by a single metastable state and the occurences to it determine the loss of memory.

 \subsubsection{Control of the thermalization}
 
  We state the following result: 
  
  \begin{lemma}\label{picg} 
   For  $t\in [R,2 R] $, any $z\in B_\a$  and any $y\in A$ let
\be{gn}
g^z_{t}(y):= \mathbf{P}(X^{z}_t=y, \t^{z}_{G}>t)-\pi_{A}(y),
\ee
and define
\be{l2}
g:=\sup_{t\in [R,2 R] }\sup_{z\in B_\a}\frac{1}{2}\sum_{y\in A}|g^z_{t}(y)|. 
\ee
The following inequality holds: 
\be{g}
g\le d + f + f^{1 - \a} +  2f^{\a}.
\ee
   \end{lemma}

 \textit{Proof} 
  
We first note that 
\be{4}
\mathbf{P}(X^{z}_t=y, \t^{z}_{G}>t)= \m^{z}_t(y)-\mathbf{P}(X^{z}_t=y, \t^{z}_{G}\le t),
\ee
and since for any $z \in B_{\a}$, $\sum_{y\in A} \mathbf{P}(X^{z}_t=y, \t^{z}_{G}\le t)\le\mathbf{P}(\t^{z}_{G}\le t)\le 2 f^\a$, Lemma \ref{picg} will 
follow by triangular inequality once we prove that 
 \be{dtvevolRext}
\sup_{t\in[R,2R]}\sup_{z \in B_{\a}} d_{TV}(\m^z_t,\p_A)\le d+f^{1-\a} + f.
\ee

     Combining Lemma \ref{deleq}, Lemma \ref{mis} and the triangular inequality, for $t \in [R,2R], z \in B_{\a}$, we obtain \eqref{dtvevolRext}. 
 
 $\qed$
 
  \subsubsection{Control of the recurrence time}

\bl{l0} Assume that $d$ and $f$ are small enough for the quantity
 \begin{equation}
  \bar c := \frac{1}{2} - \sqrt{\frac{1}{4} - (r+ 2f^{\alpha})}
 \end{equation} 
 to be well defined. Let $t\ge 2R$ and $t'\le R$; for any $x\in B_{\alpha}$, one has the inequalities:

\be{6}
1\le { \mathbf{P}(\t^{x}_{G}>t-t')\over \mathbf{P}(\t^{x}_{G}>t)}\le 1+ 2 \bar c,
\ee
\be{6pa}
1\le { \mathbf{P}(\t^{\p_A}_{G}>t-t')\over \mathbf{P}(\t^{\p_A}_{G}>t)}\le 1+ 2 \bar c.
\ee

\el
 \textit{Proof}
  
 The lower bound is trivial. As far as the upper bound is concerned, by monotonicity  it is sufficient to
prove that
\be{7}
\frac{\mathbf{P}(\t^{x}_{G}>t-R)}{\mathbf{P}(\t^{x}_{G}>t)} \leq  1+ 2\bar c.
\ee
Define
$\s:=\s^{x}(t-2R)$ (see definition (\ref{st})) and apply Markov's property to get
$$
\mathbf{P}(\t^{x}_{G}\in(t-R,t))= \mathbf{P}(\t^{x}_{G}\in(t-R,t), \s> t-R)+
\mathbf{P}(\t^{x}_{G}\in(t-R,t),\s\le t-R)
$$
$$
\leq \mathbf{P}(\t^{x}_{G} > t - 2R) \left[ \sup_{z \in \cX} \mathbf{P}(\t^{z}_{G \cup B_{\alpha}} \ge R ) + \sup_{z \in B_{\alpha}} \mathbf{P}(\t^{z}_{G} \leq 2R)\right]
$$
$$
\leq \mathbf{P}(\t^{x}_{G} > t - 2R)\left[ r + 2 f^{\alpha} \right].
$$
where we made use of the definition of $B_{\alpha}$ and of inequality \eqref{ricB}. 
Defining 
$$c := r + 2 f^{\alpha},$$ we  get the inequality:
 \be{7bis}
  \mathbf{P}(\t^{x}_{G} > t - R) \leq \mathbf{P}(\t^{x}_{G} > t ) + c \mathbf{P}(\t^{x}_{G} > t - 2R),
 \ee
from which we deduce: 
 \be{recy}
  1 - c\frac{\mathbf{P}(\t^{x}_{G} > t - 2R)}{\mathbf{P}(\t^{x}_{G} > t - R)} \leq \frac{\mathbf{P}(\t^{x}_{G} > t )}{\mathbf{P}(\t^{x}_{G} > t - R)}.
 \ee

 For $i \geq 1$, define the quantity 
 \begin{equation}
  y_{i} := \frac{\mathbf{P}(\t^{x}_{G} > iR)}{\mathbf{P}(\t^{x}_{G} > (i-1)R)}. 
 \end{equation} 
  
  Equation \eqref{recy} entails that the sequence $(y_{i})_{i \geq 1}$ satisfies the recursive inequality:
 \begin{equation*}\label{recy1}
  1 - \frac{c}{y_{i-1}} \leq y_{i}. 
 \end{equation*} 

 This implies that for every $i \geq 1$, one has
 \begin{equation}\label{min}
  y_{i} \geq 1-\overline{c}. 
 \end{equation} 

 Let us show \eqref{min} by induction; if there exists $i \geq 1$ such that $y_{i} \geq 1-\overline{c}$, then 
for every $j \geq i$, one has $y_{j} \geq 1-\overline{c}$. Indeed, making use of equation \eqref{recy1}, we get:
 \begin{equation*}
  y_{i+1} \geq 1 - \frac{c}{y_{i}} \geq 1 - \frac{c}{1-\overline{c}} = 1- \overline{c}.
 \end{equation*} 

 For the base of the induction, we note that 
 \begin{equation}
  y_{1} = \mathbf{P}(\t^{x}_{G} > R) \geq 1 - \mathbf{P}(\t^{x}_{G} \leq 2R) \geq 1 - \sup_{z \in B_{\alpha}} \mathbf{P}(\t^{z}_{G} \leq 2R)
 \geq 1 - 2f^{\alpha} > 1- c > 1- \overline{c}.
 \end{equation} 

 In particular, this implies that for every $t = kR, k \geq 2$: 
 \begin{equation}
 \frac{\mathbf{P}(\t^{x}_{G}>t-R)}{\mathbf{P}(\t^{x}_{G}>t)} = \frac{1}{y_{k}} 
\leq \frac{1}{1-\overline{c}}  \leq 1+2 \bar c,
 \end{equation} 
since $\bar c<{1\over 2}$, which entails the claim of equation \eqref{7} in the case 
when $t/R$ is an integer. 

  Assume now that $t/R$ is not an integer, and define $k = \lfloor t/R \rfloor$, so that $t=kR+t_0$
 with $t_0< R$ and define
 $$
  y_{i}(t_0) := \frac{\mathbf{P}(\t^{x}_{G} > iR+t_0)}{\mathbf{P}(\t^{x}_{G} > (i-1)R+t_0)}. 
 $$
 For each $t_0<R$ we have the same recursive inequality for $y_{i}(t_0)$ and again for the base of the induction
 we have
 $$
 y_{1}(t_0) := \frac{\mathbf{P}(\t^{x}_{G} > R+t_0)}{\mathbf{P}(\t^{x}_{G} > t_0)}\ge \mathbf{P}(\t^{x}_{G} > 2R)
 $$
 so that the induction still holds;
 this concludes the proof of \eqref{7} for every $t\ge 2R$. 
 
The proof of (\ref{6pa}) follows the same way, noting that one uses condition \eqref{f1} instead of \eqref{defB} to initialize the recurrence. 

 $\qed$ 

   \subsubsection{ Proof of Theorem \ref{convdtv} for starting points $z\in B_\a$}

As before, considering the time $\s^x=\s^x(t-2R)$, the following equality holds:
$$
\mathbf{P}(X^{x}_t=y, \t^{x}_{G}>t)
$$
$$
=
\sum_{z\in B_\a}
\int_{t-2R}^{ t-R}
\mathbf{P}(X^{x}_s=z, \t^{x}_{G}>s, \;\s^x\in(s,s+ds))\mathbf{P}(X^{z}_{t-s}=y, \t^{z}_{G}>t-s)
$$
$$+\mathbf{P}(X^{x}_t=y, \t^{x}_{G}>t, \;\s^x>t-R).
$$
 
By the definition of $g^z_{t}(y)$, we get
$$
\mathbf{P}(X^{x}_t=y, \t^{x}_{G}>t)=\sum_{z\in B_\a}\int_{t-2R}^{ t-R}
\mathbf{P}(X^{x}_s=z, \t^{x}_{G}>s, \;\s^x\in(s,s+ds))\Big[ \pi_{A}(y)+g^z_{t-s}(y)\Big]
$$
$$+\mathbf{P}(X^{x}_t=y, \t^{x}_{G}>t, \;\s^x>t-R)
$$
so that
$$
\sum_{y\in A} | \tilde\mu^{x}_t(y)-\pi_{A}(y)|\le \sum_{y\in A}
\Big| { \mathbf{P}( \t^{x}_{G}>\s^{x}, \;\s^{x}\in[t-2R,t-R]) \over \mathbf{P}(\t^{x}_{G}>t)}  -1\Big|
\pi_{A}(y)$$
$$+\sum_{y \in A } \sum_{z\in B_\a}\int_{t-2R}^{ t-R}
{\mathbf{P}(X^{x}_s=z, \t^{x}_{G}>t-2R, \;\s^{x}\in(s,s+ds))\over \mathbf{P}(\t^{x}_{G}>t)}|g^z_{t-s}(y)|
$$
\be{60}
+
{\mathbf{P}(\t^{x}_{G}>t-2R, \s^{x}>t-R)  \over \mathbf{P}(\t^{x}_{G}>t)}=:I+II+III.
\ee
 
By using the Markov property, the estimates \eqref{6}, \eqref{g} and the recursion in $B_\a$ Rc(R,r),
we conclude the proof by estimating the three terms, I, II and III of the r.h.s. of (\ref{60}). Note that in what follows, we will 
make repeated use of the monotonicity property. 

The first term can be estimated by
$$I\le 2 \bar c\vee r(1+2 \bar c)\le  2 \bar c\vee 2r \le  2 \bar  c.$$
 
Indeed, recalling that $t \geq 2R$, one can apply Lemma \ref{l0} to get:
$$
{ \mathbf{P}( \t^{x}_{G}>\s^x, \;\s^x\in[t-2R,t-R]) \over \mathbf{P}(\t^{x}_{G}>t)} \le
{ \mathbf{P}( \t^{x}_{G}>t-2R) \over \mathbf{P}(\t^{x}_{G}>t)} \le 1+ 2 \bar c. 
$$
 
On the other hand, making use of the Markov property at time $t-2R$: 
$$
{ \mathbf{P}( \t^{x}_{G}>\s^x, \;\s^x\in[t-2R,t-R]) \over \mathbf{P}(\t^{x}_{G}>t)} \ge
{ \mathbf{P}( \t^{x}_{G}>t-R, \s^x\in[t-2R, t-R]) \over \mathbf{P}(\t^{x}_{G}>t)} 
$$
$$\ge
{ \mathbf{P}( \t^{x}_{G}>t-R) \over \mathbf{P}(\t^{x}_{G}>t)} -
{ \mathbf{P}( \t^{x}_{G}>t-2R, \s^x> t-R) \over \mathbf{P}(\t^{x}_{G}>t)} 
\ge 1-{ \mathbf{P}( \t^{x}_{G}>t-2R) \over \mathbf{P}(\t^{x}_{G}>t)} \sup_{z\in\cX}\mathbf{P}(\t^z_{B_\a\cup G}>R)
$$
$$
\ge 1-r(1+ 2 \bar c). 
$$
 
To deal with II, we use Lemma \ref{picg} to get 
$$
II \le g { \mathbf{P}( \t^{x}_{G}>t-2R) \over \mathbf{P}(\t^{x}_{G}>t)}\le \left(2f^{\a} + f^{1_\a} +f + d\right)(1+2 \bar c)\le 4(2f+f^{\a} +d),
$$
and similarly for III:
$$
III \le r(1+2 \bar c)\le 2\bar c. 
$$

Finally, we prove (\ref{dtvpa}); from \eqref{piB}, we obtain that for $t \geq 2R$: 
$$
d_{TV}(\tilde\mu^{\p_A}_{A,t},\pi_{A})
={1\over 2}\sum_{y \in \cX}\Big| \sum_{x \in A}\p_A(x)\Big[\tilde \m^x_{A,t}(y)-\p_A(y)  \Big]  \Big|\le 10\bar c
\p_A(B_\a)+\p_A(B_\a^c)\le 10 \bar c+f^{1-\a},
$$
and finally we can conclude by making use of Remark \ref{hatB}.  


\subsubsection{ Proof of Theorem \ref{convdtv} for  starting points $x\in A\backslash B_\a$}

 We consider $t \geq 2R$, $x \in A \setminus B_{\a}$, $y \in A$ and we  make use of Markov's property to get the equality: 
\begin{align}\label{1der}
 & \mathbf{P}(X^x_t=y, \t^x_G>t, \t^x_{G\cup B_\a}<t-2R)= \nonumber \\
& \phantom{ii} \int_{0}^{t-2R}\sum_{z\in B_\a}\mathbf{P}(X^x_s=z, \t^x_G>s, \t^x_{G\cup B_\a}\in ds)
\mathbf{P}(X^z_{t-s}=y, \t^z_G>t-s).
\end{align}

  Since for $z\in B_\a$, $\tilde\mu^{z}_{A,t}$ and $\hat\mu^{z}_{A,t}$ coincide ( see Remark \ref{hatB}), making use of
  \eqref{dtv} for starting points 
   $z \in B_{\a}$,
   which we already proved in the previous section, we get that for $s \in [0,t-2R]$, the quantity $f_{s}^z(y)$ defined by 
  \begin{equation}
    f_{s}^z(y) =  \tilde \mu^{z}_{t-s,A}(y) - \p_A(y) 
  \end{equation} 
 satisfies 
  \begin{equation}\label{2der}
  \sup_{s \leq t-2R} \sup_{z\in B_\a} \sum_{y \in A}  \left| f_{s}^z(y) \right| <\e_{1}.
  \end{equation} 
  
   On the other hand, as in \eqref{1der}, we get:
 \begin{align}\label{3der}
    & \mathbf{P}(\t^x_G>t, \t^x_{G\cup B_\a}<t-2R)= \nonumber \\
& \phantom{ii} \int_{0}^{t-2R}\sum_{z\in B_\a}\mathbf{P}(X^x_s=z, \t^x_G>s, \t^x_{G\cup B_\a}\in ds)
\mathbf{P}(\t^z_G>t-s).
 \end{align}
 
  Combining \eqref{1der} and \eqref{3der}, we can write:
 \begin{align}
&\mathbf{P}(X^x_t=y, \t^x_G>t, \t^x_{G\cup B_\a}<t-2R)-\p_A(y) \mathbf{P}( \t^x_G>t, \t^x_{G\cup B_\a}<t-2R)\\
&=\int_{0}^{t-2R}\sum_{z\in B_\a}\mathbf{P}(X^x_s=z, \t^x_G>s, \t^x_{G\cup B_\a}\in ds)\Big(\p_A(y) \mathbf{P}(\t^z_G>t-s)+ f_{s}^z(y)\mathbf{P}(\t^z_G>t-s)\Big)\\
& \phantom{iii} -\p_A(y)\int_{0}^{t-2R}\sum_{z\in B_\a}\mathbf{P}(X^x_s=z, \t^x_G>s, \t^x_{G\cup B_\a}\in ds)
\mathbf{P}(\t^z_G>t-s)\\
&=\int_{0}^{t-2R}\sum_{z\in B_\a}\mathbf{P}(X^x_s=z, \t^x_G>s, \t^x_{G\cup B_\a}\in ds) f_{s}^z(y)\mathbf{P}(\t^z_G>t-s).
\end{align}
  
   From the last equality and \eqref{3der}, we deduce that for any $t \geq 2R$, $x \in A \setminus B_{\a}$ and $y \in A$, 
  \begin{equation}
   \left|\hat \mu^{x}_{t,A}(y) - \pi_{A}(y) \right| \leq \sup_{s \leq t-2R} \sup_{z \in B_{\a}} |f_{s}^z(y)|,
  \end{equation} 
  from which \eqref{dtv} follows from considering \eqref{2der}.

%
%


\subsection{Proof of Theorem \ref{expbeh}}

We define
   \begin{equation}
    T^{*} = \mathbf{E}[\tau_{G}^{\mu^{*}_{A}}].
   \end{equation} 
   
    For $y \in A$, we consider the quantity
   \begin{equation}
    \delta(y):=2d_{TV}(\hat{\mu}_{A,2R}^{y},\mu_{A}^{*}),               
             \end{equation} 
  and recalling \eqref{dtvqs}, we have
  $$
 \delta := \sup_{y\in A}\d(y)\le \e_1+\e_2.
  $$

   We first show that the recurrence time is asymptotically negligible with respect to 
   $T^{*}$.

    \begin{lemma}\label{rat}
    There exists a constant $C > 0$ such that the following inequality holds:
     \begin{equation}
      \frac{R}{T^{*}} \leq C(f + \delta).
     \end{equation} 
    \end{lemma}
    
     \textit{Proof}
     
      By Proposition \ref{exppro}, we have:
      \begin{equation}
       \mathbf{P}[\tau_{G}^{\mu^{*}_{A}} > 2R] = e^{-2R/T^{*}}.
      \end{equation} 
     
      On the other hand, by making use of \eqref{dtvqs} and of $E(R,f)$:
       \begin{align}
         \mathbf{P}[\tau_{G}^{\mu^{*}_{A}} > 2R] & = \sum_{z \in A} \mu^{*}_{A}(z) \mathbf{P}[\tau_{G}^{z} > 2R] 
         = \mathbf{P}[\tau_{G}^{\pi_{A}} > 2R] + \sum_{z \in A} \left[ \mu^{*}_{A}(z) - \pi_{A}(z)\right] \mathbf{P}[\tau_{G}^{z} > 2R]  \\
          & \geq \mathbf{P}[\tau_{G}^{\pi_{A}} > 2R] - d_{TV}( \pi_{A},\mu^{*}_{A}) \geq 1 - f - \delta,
       \end{align}
        which proves the claim.
      $\qed$
      
       \begin{lemma}\label{recT}
       The following holds for any $k \geq 0$:
\begin{equation}
 \sup_{y\in A}\mathbf{P}(\tau_{G}^{y}>2kR)\le e^{-k\frac{2R}{T^{*}}}\left(1+O\left(\delta+r+\frac{R}{T^{*}}\right)\right).
\end{equation} 
     \end{lemma}

  \textit{Proof}
  Making use of Markov's property at time $2R$, we have for any $t \geq 2R$ and $y \in A$
  \begin{align}
     \mathbf{P}(\tau_{G}^{y}>t)& =  \mathbf{P}(\tau_{G}^{y}>2R, \t^y_{B_\a\cup G}\le R)
    \sum_{z\in A}\hat\m_{A,2R}^{y}(z)\mathbf{P}(\tau_{G}^{z}>t-2R)\\
  &+ 
    \sum_{z\in A} \mathbf{P}(\tau_{G}^{y}>2R, \t^y_{B_\a\cup G}>R,  X^y_{2R}=z)    \mathbf{P}(\tau_{G}^{z}>t-2R). 
 \end{align}
     
   By using  Proposition \ref{exppro} and the hypothesis $Rc(R,r)$, we have the following estimate:
     \begin{align*}
     & \mathbf{P}(\tau_{G}^{y}>t) \le \mathbf{P}(\tau_{G}^{y}>2R, \t^y_{B_\a\cup G}\le R)\left[e^{-\frac{t-2R}{T^{*}}}
    +\sum_{z\in A}\left(\hat{\mu}_{A,2R}^{y}(z)-\mu^{*}_{A}(z)\right)\mathbf{P}(\tau_{G}^{z}>t-2R)\right] \\ 
    & \phantom{iiiiiiiiiiiiiiiiiii} + r \sup_{z\in A}\mathbf{P}(\tau_{G}^{z}>t-2R).
     \end{align*} 
      
      In particular, for $t = 2(k+1)R$, we get:
       \begin{equation}\label{eq:2}
   \Big|   
   \mathbf{P}(\t^y_G>(k+1)2R)-\mathbf{P}(\tau_{G}^{y}>2R, \t^y_{B_\a\cup G}\le R)e^{-\frac{t-2R}{T^{*}}} \Big|
   \le (\d+r) \sup_{z\in A}  \mathbf{P}(\t^z_G>k2R)
  \end{equation} 
  and 
   \begin{equation}\label{eq:3}
    \sup_{y\in A}\mathbf{P}(\tau_{G}^{y}>2(k+1)R) \leq e^{-k\frac{2R}{T^{*}}}+(\delta+r)\sup_{z\in A}\mathbf{P}(\tau_{G}^{z}>2kR).
   \end{equation}

     From \eqref{eq:3}, we prove by recurrence that 
    \begin{equation}\label{eq:r4}
\sup_{y\in A}\mathbf{P}(\tau_{G}^{y}>2kR)\le (\delta+r)^{k}
+\frac{e^{-k\frac{2R}{T^{*}}}}{e^{-\frac{2R}{T^{*}}}-(\delta+r)}.
   \end{equation} 
   
    Indeed, by Lemma \ref{rat}, as soon as $R/T^{*}, \delta$ and $r$ are small enough, \eqref{eq:r4} holds for $k=0$; assume
    \eqref{eq:r4} to be true for a given $k \geq 0$. Using \eqref{eq:3}, we get that
    \begin{align}
     \sup_{y\in A}\mathbf{P}(\tau_{G}^{y}>2(k+1)R) & \le e^{-(k+1)\frac{2R}{T^{*}}}
     + (\delta+r) \left( (\delta+r)^{k}+\frac{e^{-k\frac{2R}{T^{*}}}}{e^{-\frac{2R}{T^{*}}}-(\delta+r)} \right) \\
   & \leq (\delta+r)^{k+1}+e^{-k\frac{2R}{T^{*}}}
   \left(\frac{e^{-\frac{2R}{T^{*}}}-(\delta+r)+(\delta+r)}{e^{-\frac{2R}{T^{*}}}-(\delta+r)}\right),
    \end{align} 
    which closes the recursion. 
    
     Now we note that \eqref{eq:r4} implies Lemma \ref{recT}. 
   
   $\qed$
   
%
%
%
     
      \textit{Proof of Theorem \ref{expbeh}}
     
   For $x \in A \setminus B_{\a}$, \eqref{fortisB} directly follows from \eqref{eq:2} and from \eqref{eq:r4}. Note that one gets 
    from a statement on a time scale $2R$ to a statement on a generic time scale $t \geq 0$ in a standard way and by using Lemma \ref{rat}. 
    
     To get the statement \eqref{fortisl}, 
    we note that for $x \in B_{\a}$, we have:
     \begin{equation}
     \left| \mathbf{P}[\tau_{G}^{x} > 2R, \tau_{B_{\a} \cup G}^{x} \leq R] - 1 \right| \leq r + f^{\a}, 
     \end{equation} 
   and hence we are done. 
   
    Finally, in the same way, making use of the slow escape property $E(R,f)$, we get that
\begin{equation}\label{misa}
 \mathbf{P}(\tau_{G}^{\pi_{A}}>t)  =  e^{-\frac{t}{T^{*}}}\left(1+O\left(\delta+\frac{R}{T^{*}}+f\right)\right).
\end{equation} 

 $\qed$

 {\bf Acknowledgments:}The authors thank A.Asselah, P. Caputo, A. Faggionato, A. Gaudilli\`ere, F. Martinelli and E.Olivieri for 
 fruitful discussions. The work of JS was 
 supported by NWO STAR grant: ``Metastable and cut-off behavior of stochastic processes''. The work of ES
  was partially supported by PRIN projects “Campi aleatori, percolazione ed evoluzione stocastica di sistemi
 con molte componenti”, “Dal microscopico al macroscopico: analisi di strutture comp- lesse e applicazioni”. The work
 of FM was supported by Universit\`{a} di "Roma Tre". FM and ES thank Eurandom, the department of 
 Mathematics and Computer Science at Technical University Eindhoven and the department of Mathematics at Utrecht 
 University for the kind hospitality. RF and JS thank University of "Roma Tre" for the kind hospitality.

  \bibliographystyle{abbrv}

  \bibliography{AAP14Bibli}

\end{document}